\documentstyle[12pt,amssymbols]{article}

\newfont{\g}{eufm10 scaled\magstep1}
\newfont{\gs}{eufm7}

\newcommand{\gtg}{\mbox{\g g}}
\newcommand{\Sym}{\mbox{\g S}}
\newcommand{\sym}{\mbox{\gs S}}
\newcommand{\gl}{\mbox{\g gl}}

\newfont{\bg}{cmr10 scaled\magstep4}
\newcommand{\bigzerol}{\smash{\hbox{\bg 0}}}
\newcommand{\bigzerou}{\smash{\hbox{\lower1.7ex\hbox{\bg 0}}}}
\newcommand{\bigast}{\smash{\hbox{\bg *}}}

\def\N{\Bbb N}

\def\Z{\Bbb Z}
\def\Q{\Bbb Q}
\def\M{{\cal M}}
\def\I{{\cal I}}
\def\Cl{{\cal A}}
\def\F{{\cal F}}
\def\<{\langle}
\def\>{\rangle}
\def\A{A^{(1)}_{r-1}}

\def\a{\alpha}
\def\b{\beta}

\def\e{\epsilon}
\def\l{\lambda}
\def\L{\Lambda}
\def\fr{\frac}
\def\s{\sum}

\def\Sa{S^{(r)}_\a}
\def\barSa{\bar{S}^{(r)}_\a}
\def\Sb{S^{(r)}_\b}

\def\varprojlim{\displaystyle{\lim_{\longleftarrow}}}
\def\sgn{\mbox{\rm sgn}}
\def\lra{\longrightarrow}
\def\vac{\mbox{\o}}
\def\uin{\mbox{
{\footnotesize
\unitlength 10pt
\begin{picture}(1,1)
\put(0,0.3){\makebox(0,0){$\cup$}}
\put(0,0){\line(0,1){0.61}}
\end{picture}}}
}

\newcounter{sect}
\def\thesect{\arabic{sect}}
\newtheorem{definition}{Definition}[sect]
\newtheorem{theorem}[definition]{Theorem}
\newtheorem{lemma}[definition]{Lemma}
\newtheorem{proposition}[definition]{Proposition}
\newtheorem{corollary}[definition]{Corollary}

\begin{document}

\thispagestyle{empty}
\vspace{1cm}
\begin{large}
\begin{center}
{\bf REDUCED SCHUR FUNCTIONS \\
AND \\
THE LITTLEWOOD-RICHARDSON COEFFICIENTS}
\end{center}
\end{large}

\vspace{1.0cm}

\begin{center}
Susumu ARIKI

{\it Division of Mathematics,Tokyo University of Mercantile Marine,} \\
{\it Etchujima 2-1-6, Koto-ku, Tokyo 135, Japan} \\

Tatsuhiro NAKAJIMA

{\it Department of Physics, Tokyo Metropolitan University}\\
{\it Minami-Ohsawa 1-1, Hachioji-shi, Tokyo 192-03, Japan}\\
and

Hiro-Fumi YAMADA

{\it Department of Mathematics, Tokyo Metropolitan University}\\
{\it Minami-Ohsawa 1-1, Hachioji-shi, Tokyo 192-03, Japan}

\end{center}

\newpage
\setcounter{sect}{0}
\begin{flushleft}
{\large\bf \S \thesect \hspace{5mm} Introduction}
\end{flushleft}
\vspace{3mm}

The present paper deals with a formula satisfied by ``$r$-reduced'' Schur
functions.

Schur functions originally appear as irreducible characters of general
linear group over the complex number field.
In this paper they are considered as weighted homogeneous polynomials
with respect to the power sum symmetric functions.
More precisely, for a Young diagram $\l$ of size $n$, the Schur function
indexed by $\l$ reads
$$
S_{\l}(t)=\sum_{\nu_1+2\nu_2+\cdots=n} \chi^{\l}(\nu)
          \frac{t_1^{\nu_1}t_2^{\nu_2}\cdots}{\nu_1!\nu_2!\cdots}, $$
where $\chi^{\l}(\nu)$ is the character value of the irreducible
representation $S^{\l}$ of the group algebra $\Q\Sym_n$, evaluated at the
conjugacy class of the cycle type
$\nu=(1^{\nu_1}2^{\nu_2}\cdots n^{\nu_n})$.
Setting $t_{jr}=0$ for $j=1,2,\ldots$ in $S_{\l}(t)$, we have the
$r$-reduced Schur function $S_{\l}^{(r)}(t)$.
The set of all $r$-reduced Schur functions spans the polynomial ring
$P^{(r)}=\Q[t_j ; j\not\equiv 0 (\mbox{mod }r) ]$.
We show that a good choice of basis elements leads to an explicit
description of all other $r$-reduced Schur functions involving the
Littlewood-Richardson coefficients.

The formula has not only a purely combinatorial meaning, but also
nice implications in two different fields.
One is about the basic representation of the affine Lie algebra $\A$.
We show that the basis in the main theorem gives in turn a weight basis of
the basic $\A$-module realized in $P^{(r)}$.
The other is about modular representations of the symmetric group.
Our explicit formula implies that the determination of the decomposition
matrices reduces to that for the basic set we give in this paper.

The paper is organized as follows.
In Section 1 we introduce generalized Maya diagrams and associated
$r$-reduced Schur functions.
Section 2 is for a combinatorics of Young diagrams.
Section 3 is devoted to the main theorem.
In Section 4 we describe weight vectors of the basic $\A$-module.
In Section 5 the formula is translated into that in the modular
representation theory.

The authors are grateful to A. Gyoja for some comments.

\stepcounter{sect}
\begin{flushleft}
{\large\bf \S \thesect \hspace{5mm}
Generalized Maya Diagrams and Schur Functions}
\end{flushleft}
\vspace{3mm}

In this section we will summarize known results on Schur functions.

Let $\L_n=\Q[x_1,x_2,\ldots,x_n]^{\sym_n}$ be the ring of
symmetric polynomials of $n$ variables.
For $m<n$ we have a surjective homomorphism $P_{nm}: \L_n \lra \L_m$
by setting $x_{m+1}=x_{m+2}=\cdots=x_n=0$.
In this way we obtain an inverse system
and denote $\L=\varprojlim \, \L_n$.
The power sum symmetric function $p_j(x)$ is by definition
$p_j(x)=\varprojlim \, p_j(x_1,\ldots,x_n)$,
where $p_j(x_1,\ldots,x_n)=x_1^j+\cdots+x_n^j$.
It is known that $\L$ is generated by ${p_j(x)}\,\,(j=1,2,\ldots)$
and these elements are algebraically independent.

Let $\a=(\a_1,\a_2,\ldots)$ be a semi-infinite sequence
with integer components satisfying $\a_{i+1}=\a_i-1$ for
sufficiently large $i$.
We shall introduce an equivalence relation on the
set of such sequences as follows:
$\a\sim \b$ if $\a_i-\b_i$ is independent of $i$.
The quotient by the equivalence relation is denoted by
$$
\M=\left\{\a=(\a_1,\a_2,\ldots)|
\mbox{$\a_{i+1}=\a_i-1$ for all $i\gg 0$} \right\}
\left/ \sim \right.
$$
and an element of $\M$ is called a generalized Maya diagram.
We let $\M_d$ be a subset of $\M$ consisting of strictly decreasing
sequences.

We now define Schur functions indexed by generalized Maya diagrams.
For each $\a\in \M$, we put
$$
S_{\a}(x)=\varprojlim S_{\a}(x_1,\ldots,x_n),
$$
where
$$
S_{\a}(x_1,\ldots,x_n)=
\fr{\det\left(x_i^{\a_j-\a_n}\right)_{1\leq i,j\leq n}}
{\det\left(x_i^{-j+n}\right)_{1\leq i,j\leq n}}
$$
for a sufficiently large $n$. We can easily see that
$S_{\a}(x_1,\ldots,x_n)$ are well-defined and that
$P_{nm}(S_{\a}(x_1,\ldots,x_n))=S_{\a}(x_1,\ldots,x_m)$ for $m<n$.

Let $\Sym_n$ be the symmetric group of degree $n$ and put
$$
\Sym_{\infty}=\bigcup_{n\geq 1}\Sym_n.
$$
By permuting the first $n$ components, $\Sym_n$ acts on $\M$.
If there exists a unique element $\sigma\in \Sym_{\infty}$ such that
$\sigma\a\in \M_d$, the sign of $\sigma$ is referred to as the
sign of $\a$ and denoted by $\sgn(\a)$.
Otherwise we put $\sgn(\a)=0$.

For $\a\in \M$ we define a sequence $\l(\a)=(\l_1,\l_2,\ldots)$
by putting $\l_j=\a_j-\a_n+(j-n)$ for a sufficiently large $n$.
The size $|\l(\a)|$ is defined by $\sum \l_j$.
We will write $|\a|$ instead of $|\l(\a)|$ for $\a\in \M$.
Note that in case of $\a\in \M_d$, $\l(\a)$ is a Young diagram.
We often identify an element of $\M_d$ with a Young diagram.
The transposed Young diagram of $\l(\a)$ is denoted by $\l(\a)'$
and the corresponding element in $\M_d$ is denoted by $\a'$, i.e.,
$\l(\a)'=\l(\a')$.
For $\a, \b_1,\ldots,\b_m\in \M$, the Littlewood-Richardson
coefficients are defined by
$$
C_{\b_1\cdots \b_m}^{\a}=\left\{
\begin{array}{ccl}
LR_{\l(\b_1)\cdots \l(\b_m)}^{\l(\a)} & &\mbox{if $\exists
\sigma$, $\tau_j\in \Sym_{\infty}$ ($j=1,\ldots,m$) such that
$\sigma \a,\,\tau_j \b_j\in \M_d$} \\
0& & \mbox{otherwise},
\end{array}
\right.
$$
where
$$
S_{\l_1}(x)\cdots S_{\l_m}(x)=
             \sum_{\mu} LR_{\l_1\cdots \l_m}^{\mu}S_{\mu}(x).
$$

The following properties of the Schur functions are well known \cite{M}.
\begin{proposition}
\begin{enumerate}
\item $S_{\sigma\a}(x)=\sgn(\sigma)S_{\a}(x)$
for any $\sigma\in \Sym_{\infty}$.
\item For any $j\geq 1$ and $\a\in \M$, the sequence $\a+l\e_i$
is in $\M$ and
$$
S_{\a}(x)p_j(x)=\s_{i\geq 1} S_{\a+j\e_i}(x),
$$
where
$\displaystyle{\e_i=(0,\ldots,0,\mathop{\breve 1}^i,0,\ldots)}$.
\item The set $\left\{S_{\a}(x) | \a\in \M_d \right\}$
gives a basis of $\L$.
\item For $\a\in \M_d$,
$$
S_{\a}(x,y)=
\sum_{\b,\gamma\in \M_d} C_{\b \gamma}^{\a} S_{\b}(x) S_{\gamma}(y).
$$
\item Let $\omega: \L \lra \L$ be an algebra automorphism
defined by $\omega(p_j(x))=(-1)^{j-1}p_j(x)$. Then
$$
\omega\left(S_{\a}(x)\right)=S_{\a'}(x) \,\,\,(\a\in \M_d).
$$
\end{enumerate}
\end{proposition}

Note that the summation in (2) is a finite sum.
If $\a\not\in M_d$, then we define $S_{\a'}(x)$ by
$S_{\a'}(x)=\omega(S_{\a}(x))$. This definition is compatible with
above (1) and (5).

Let $\rho^{(r)}:\L=\Q[p_1,p_2,\ldots] \lra
P^{(r)}=\mbox{$\Q[t_j ; j\not\equiv 0$ (mod $r$)}]$
be the surjective homomorphism of algebras defined by
$$
\left\{
\begin{array}{cccl}
p_j(x) & \longmapsto &  jt_j & \mbox{if $j\not\equiv 0$ (mod $r$)} \\
& & & \\
p_j(x) & \longmapsto &  0    & \mbox{otherwise}.
\end{array}
\right.
$$
We shall write $\Sa(t)=\rho^{(r)}\left(S_{\a}(x)\right)$ and call it the
$r$-reduced Schur function indexed by $\a$.
When $\a\in \M_d$  we may write $S_{\l(\a)}^{(r)}(t)$ instead of $\Sa(t)$.
The following are easily checked by using Proposition 1.1(1) and (2).
\begin{lemma}
(1)\hspace{12pt} $\displaystyle{S_{\sigma\a}^{(r)}(t)=\sgn(\sigma)\Sa(t)}.$
\hspace{1cm}(2)\hspace{12pt}
$\displaystyle{\sum_{i\geq 1} S_{\a+jr\e_i}^{(r)}(t)=0}$ $(j\geq 1)$.
\end{lemma}

\stepcounter{sect}
\begin{flushleft}
{\large\bf \S \thesect \hspace{5mm}
Cores and Quotients}
\end{flushleft}
\vspace{3mm}

Let $\a\in \M$ be a generalized Maya diagram.
Define, for $0\leq k\leq r-1$,
$$
\a^{(k)}=\{\a_i | \a_i\equiv k (\mbox{mod }r)\},
\a_+^{(k)}=\left\{\a_i\in \a^{(k)} | \a_i \geq 0\right\},
\a_-^{(k)}=\left\{\a_i\in \a^{(k)} | \a_i < 0\right\}.
$$
We fix a sequence in the equivalence class $\a$ satisfying the
following condition: \linebreak
(1) $\displaystyle{\bigcup_k \a_-^{(k)}=\{-1,-2,\ldots\}}$,  \hspace{3mm}
(2) cardinality $\sum_k n_k$ is divisible by $r$, where $n_k=\#\a_+^{(k)}$.

We can form new elements $\a[k]\in \M$ and $\a_c\in \M_d$ as follows.
The sequence $\a[k]$ is obtained by
$$
\a[k]=\left(\frac{\a_{i_1}^{(k)}-k}{r},\frac{\a_{i_2}^{(k)}-k}{r},
\ldots \right),
$$
where $(\a_{i_1}^{(k)},\a_{i_2}^{(k)},\ldots)$
is the subsequence consisting of elements of $\a^{(k)}$.
We have $\a_c$ if we replace $(\a_{i_1}^{(k)},\a_{i_2}^{(k)},\ldots)$ by
$(n_k r+k,(n_k-1)r+k,\ldots,k,-r+k,\ldots)$.

In the case $\a\in \M_d$, $\l(\a_c)$ and $(\l(\a[0]),\ldots,\l(\a[r-1]))$
is called the $r$-core and the $r$-quotient of $\l(\a)$,
respectively  \cite{O}.
The $r$-core is obtained from $\l$ by
removing $r$-hooks successively as many as possible.

Since non-negative entries of $(\a_1^{(0)},\ldots,\a_1^{(r-1)},
\a_2^{(0)},\ldots,\a_2^{(r-1)},\ldots)$ is a permutation of that of
$(\a_1,\a_2,\ldots)$, we define the $r$-sign $\delta_r(\a)$ of
$\a\in \M_d$ to be the sign of this permutation.
For an arbitrary $\a\in \M$, we define
$$
\delta_r(\a)=\left\{
\begin{array}{cl}
\delta_r(\sigma\a) &
\mbox{if $\exists \sigma\in \Sym_{\infty}$ such that $\sigma\a\in \M_d$} \\
0 & \mbox{otherwise}.
\end{array} \right.
$$

If $r=2$ and $\a\in \M_d$, then $\l(\a_c)$ is a staircase Young diagram
and $\delta_2(\a)$ is simply described as $(-1)^q$,
where $q$ is the number of column 2-hooks to be removed in
the procedure of completing the 2-core of $\l(\a)$.

If $\a\in \M$ corresponds to the $(r+1)$-tuple $(\a_c;\a[0],\ldots,\a[r-1])$
we may also write \linebreak
$S^{(r)}_{(\a_c;\a[0],\ldots,\a[r-1])}(t)$ instead of $\Sa(t)$.
We have
$$
S^{(r)}_{(\a_c;\sigma_0\a[0],\ldots,\sigma_{r-1}\a[r-1])}(t)
=\prod_{k=0}^{r-1} \sgn(\sigma_k) S^{(r)}_{(\a_c;\a[0],\ldots,\a[r-1])}(t)
\mbox{\hspace{5mm} for $\sigma_k\in \Sym_{\infty}$}.
$$

\stepcounter{sect}
\begin{flushleft}
{\large\bf \S \thesect \hspace{5mm} Main Theorem}
\end{flushleft}
\vspace{3mm}

We are now ready to state our main theorem.
\begin{theorem}
For any Young diagram $\l$ we have
$$
S^{(r)}_{\l}(t)=(-1)^{|\l[0]|} \delta_r(\l)
\! \! \sum_{\mu,\nu_1,\ldots,\nu_{r-1}} \! \! \!
LR^{\l[0]'}_{\nu_1\cdots \nu_{r-1}} LR^{\mu[1]}_{\nu_1 \l[1]} \cdots
LR^{\mu[r-1]}_{\nu_{r-1} \l[r-1]} \delta_r(\mu)
S^{(r)}_{\mu}(t),
$$
where summation runs over Young diagrams $\mu$ and
$\nu_1,\ldots,\nu_{r-1}$ such that $|\mu|=|\l|$, $\mu[0]=\vac$ and
the core of $\mu$ coincides with
that of $\l$.
\end{theorem}

Let $\theta$ denote the element of $\M_d$ corresponding to
the empty Young diagram $\vac$.
We shall show that, for any $\a\in \M$,
$$
S^{(r)}_{\a}(t) \\
=(-1)^{|\a[0]|} \sgn(\a) \delta_r(\a)
\! \! \! \sum_{\beta,\gamma_1,\ldots,\gamma_{r-1}} \! \! \! \! \!
C^{{\a[0]}'}_{\gamma_1\cdots\gamma_{r-1}} C^{\b[1]}_{\gamma_1 \a[1]}
\cdots C^{\b[r-1]}_{\gamma_{r-1}\a[r-1]}
\delta_r(\b) \Sb(t),
$$
where summation runs over $\beta, \gamma_1,\ldots,\gamma_{r-1}\in \M_d$
such that $|\beta|=|\a|$,
$\beta[0]=\theta$ and $\beta_c=\a_c$.
Let $F_{\a}(t)$ be the right hand side.
We focus on identities satisfied by $F_{\a}(t)$ and $\Sa(t)$.
We first need a lemma.
\begin{lemma}
Let $U$ be the vector space
$$
U:=\bigoplus_{\a\in \M}\Q u_{\a}\left/\sum_{\a\in \M, \sigma\in \Sym_{\infty}}
\! \! \! \! \! \Q(u_{\sigma\a}-\sgn(\sigma)u_{\a})\right..
$$
Then we have the following:
\begin{enumerate}
\item $U$ has a basis $\left\{u_{\a}| \a\in M_d \right\}$.
\item If we set
$$
U_1:=U\left/\sum_{j\geq 1} \Q \left(\sum_{i\geq 1} u_{\a+j\e_i}\right),\right.
$$
then $U_1 \simeq \Q$.
\item If we set
$$
U_r:=U\left/\sum_{j\geq 1} \Q \left(\sum_{i\geq 1} u_{\a+jr\e_i}\right),
\right.
$$
then there is a linear isomorphism
$$
\begin{array}{ccc}
U_r & \widetilde{\lra} & P^{(r)} \\
\uin  &                & \uin \\
u_{\a} & \longmapsto & \Sa(t).
\end{array}
$$
\end{enumerate}
\end{lemma}
\noindent
{\bf Proof.}  (1) is obvious.

(2) We have a canonical linear surjection from $\L$ to $U_1$
which maps $S_{\a}$ to $u_{\a}$.
Then the kernel of this surjection coincides with the maximal ideal
$\I=(p_1(x),p_2(x),\ldots)$ of $\L$
because of the relations in Proposition 1.1(2) and
$\sum_{j\geq 1} u_{\a+j\e_i}=0$ in $U_1$.
Since the algebra $\L/\I$ is isomorphic to $\Q$, we have
$U_1 \simeq \L/\I \simeq \Q$.

(3) Consider a linear surjection
$T: \L \lra U_r$ defined by $T(S_{\a}(x))=u_{\a}$.
Put $\I^{(r)}=(p_r(x),p_{2r}(x),\ldots)$.
Then, by using Proposition 1.1(2), we have $T(\I^{(r)})=0$
and a surjection $\bar{T}: \L/\I^{(r)} \lra U_r$.
On the other hand we can define a linear surjection
$S:U_r \lra P^{(r)}$ by $S(u_{\a})=\Sa(t)$.
The composition $S\circ \bar{T}$ gives a linear isomorphism
from $\L/\I^{(r)}$ to $P^{(r)}$.
Hence we have that $S$ is a linear isomorphism as desired.
$\blacksquare$

This lemma leads to the following proposition.
\begin{proposition}
The set $\left\{\Sa(t) \mid \a\in \M_d,\, \a[0]=\theta \right\}$
gives a basis of $P^{(r)}$.
\end{proposition}
\noindent
{\bf Proof. }
We first show that these elements span $P^{(r)}$.
Introduce a filtration $\{V_n\}_{n\geq 0}$ of $P^{(r)}$ by
$$V_n=\sum_{\a;|\a[0]|\leq n} \Q \Sa(t),$$
and put
$$\bar{V}=\bigoplus_{n\geq 0} V_n/V_{n-1},$$
where $V_{-1}=\{ 0 \}$.
The set $\{ \barSa(t) | \a\in \M_d,\,|\a[0]|=n \}$ spans the
direct summand $V_n/V_{n-1}$.
Here $\barSa(t)$ stands for the modulo class represented by $\Sa(t)$.
The equation in Lemma 1.2(2) reads
$$
0=\sum_{i\geq 1} S_{\a+jr\e_i}^{(r)}(t)
=\sum_{i\geq 1} S^{(r)}_{(\a_c;\a[0]+j\e_i,\a[1],\ldots,\a[r-1])}(t)
+\sum_{k=1}^{r-1} \sum_{i\geq 1}
             S^{(r)}_{(\a_c;\a[0],\ldots,\a[k]+j\e_i,\ldots,\a[r-1])}(t).
$$
If $|\a[0]+j\e_i|=n$, then
$\bar{S}^{(r)}_{(\a_c;\a[0],\ldots,\a[k]+j\e_i,\ldots,\a[r-1])}(t)=0$.
Hence
$$
\sum_{i\geq 1} \bar{S}^{(r)}_{(\a_c;\a[0]+j\e_i,\a[1],\ldots,\a[r-1])}(t)=0.
$$
If we set
$\bar{W}=\sum_{\b\in \M} \bar{S}^{(r)}_{(\a_c;\b,\a[1],\ldots,\a[r-1])}(t)$,
then it turns out to be
$\bar{W}=\Q \bar{S}_{(\a_c;\theta,\a[1],\ldots,\a[r-1])}(t)$
by applying Lemma 3.2(2).
Thus we have the first part of the proof.

Next we will show the linear independence.
We make $P^{(r)}$ into a graded algebra by putting $\deg t_j=j$.
The dimension of the homogenous component of degree $n$ is
$p^{(r)}(n)$, the number of partitions of $n$ into positive integers
not divisible by $r$.
We denote by $d^{(r)}(n)$ the cardinality of the set
$\left\{\Sa(t) \mid \a\in \M_d,\, \a[0]=\theta, \deg \Sa(t)=n \right\}$.
By the one-to-one correspondence between the set of all Young diagrams
and the set of
$(r+1)$-tuples of $r$-cores and $r$-quotients, we have
$$
\sum_n d^{(r)}(n) q^n=\frac{\phi(q^r)}{\phi(q)}=\sum_n p^{(r)}(n) q^n,
$$
where $\phi(q)=\prod_{j=1}^{\infty}(1-q^j)$.
This proves the linear independence.
$\blacksquare$

We now look at identities satisfied by
$$
F_{\a}(t)=
(-1)^{|\a[0]|} \sgn(\a) \delta_r(\a)
\! \! \! \sum_{\beta,\gamma_1,\ldots,\gamma_{r-1}} \! \! \! \! \!
C^{{\a[0]}'}_{\gamma_1\cdots\gamma_{r-1}} C^{\b[1]}_{\gamma_1 \a[1]}
\cdots C^{\b[r-1]}_{\gamma_{r-1}\a[r-1]}
\delta_r(\b) \Sb(t).
$$

\begin{proposition} \hspace{5mm}
$\displaystyle{\sum_{i\geq 1} F_{\a+jr\e_i}(t)=0}$ $(j\geq 1).$
\end{proposition}
\noindent
{\bf Proof.}
We shall prove, for any $\b$,
$$
\sum_{i\geq 1}(-1)^{|(\a+jr\e_i)[0]|} \sgn(\a+jr\e_i) \delta_r(\a+jr\e_i)
\! \! \! \! \sum_{\gamma_1,\ldots,\gamma_{r-1}} \! \! \! \! \!
C^{(\a+jr\e_i)[0]'}_{\gamma_1\cdots\gamma_{r-1}}
C^{\b[1]}_{\gamma_1 (\a+jr\e_i)[1]}\cdots
C^{\b[r-1]}_{\gamma_{r-1}(\a+jr\e_i)[r-1]}=0.
$$
Since
$$
\frac{\sgn(\a+jr\e_i)\delta_r(\a+jr\e_i)}{\sgn(\a)\delta_r(\a)}=
\frac{\sgn(\a[k]+j\e_i)}{\sgn(\a[k])}
$$
for some $k$, it is enough to show that
\begin{eqnarray*}
\lefteqn{
\sum_{i\geq 1}\left\{(-1)^j \sgn(\a[0]+j\e_i)
C^{(\a[0]+j\e_i)'}_{\gamma_1\cdots\gamma_{r-1}} C^{\b[1]}_{\gamma_1 \a[1]}
\right.\cdots C^{\b[r-1]}_{\gamma_{r-1}\a[r-1]}} \\
&+&\sum_{k=1}^{r-1} \sgn(\a[k]+j\e_i)
\left.C^{{\a[0]}'}_{\gamma_1\cdots\gamma_{r-1}}
C^{\b[1]}_{\gamma_1 \a[1]} \cdots C^{\b[k]}_{\gamma_k (\a[k]+j\e_i)}
\cdots C^{\b[r-1]}_{\gamma_{r-1}\a[r-1]} \right\}=0.
\end{eqnarray*}
Taking a summation over $\b[1],\ldots,\b[r-1]$ after multiplying by
$\displaystyle{\prod_{k=1}^{r-1}} S_{\b[k]}(x^{(k)})$,
it reduces to
\begin{eqnarray*}
\lefteqn{
\sum_{i\geq 1} (-1)^j S_{(\a[0]+j\e_i)\prime}(x^{(1)},\ldots,x^{(r-1)})
\prod_{k=1}^{r-1} S_{\a[k]}(x^{(k)})} \\
&+&
\sum_{k=1}^{r-1} \sum_{i\geq 1} S_{(\a[0])\prime}(x^{(1)},\ldots,x^{(r-1)})
S_{\a[1]}(x^{(1)})\cdots S_{\a[k]+j\e_i}(x^{(k)})\cdots
S_{\a[r-1]}(x^{(r-1)})=0.
\end{eqnarray*}
Applying Proposition 1.1(2) and (5) we see
\begin{eqnarray*}
\lefteqn{(-1)^{2j-1}S_{\a[0]\prime}(x^{(1)},\ldots,x^{(r-1)})
p_j(x^{(1)},\ldots,x^{(r-1)})
\prod_{k=1}^{r-1} S_{\a[k]}(x^{(k)})} \\
&+&\sum_{k=1}^{r-1}
S_{\a[0]\prime}(x^{(1)},\ldots,x^{(r-1)})
p_j(x^{(k)}) \prod_{l=1}^{r-1} S_{\a[l]}(x^{(l)})=0.
\mbox{\hspace{12pt}}\blacksquare
\end{eqnarray*}

Finally we complete the proof of the main theorem.
By Lemma 3.2(3) we have a linear mapping $P^{(r)} \lra P^{(r)}$,
which maps $\Sa(t)$ to $F_{\a}(t)$.
It turns out to be the identity since $\Sa(t)$ equals $F_{\a}(t)$
for the basis given in Proposition 3.3.

\stepcounter{sect}
\begin{flushleft}
{\large\bf \S \thesect \hspace{5mm} Weight Vectors in the Basic $\A$-Module}
\end{flushleft}
\vspace{3mm}

We can explicitly obtain a basis of each weight space of the basic
$\A$-module \cite{DJKM,K,KKLW} by using the results in Section 3.

Consider the affine Lie algebra  of type $\A$.
The Cartan subalgebra is
$\oplus_{i=0}^{r-1} \Q \a_i^{\vee} \oplus \Q d_0$, where
$\<d_0,\alpha_i\>=1$ for $0\leq i\leq r-1$.
The basic $\A$-module is the simple highest weight module with highest
weight $\L_0$ where $\<\a_i^{\vee}, \L_0\>=\delta_{i0}$ and
$\<d_0, \L_0\>=0$.

Let $W$ be a vector space over $\Q$ with a basis
$\{\psi(i), \psi^{*}(i) \mid i\in \Z \}$.
The Clifford algebra $\Cl$ is an associative algebra generated by $W$
with respect to the following relations:
\begin{eqnarray*}
&&\{\psi^*(i),\psi(j)\}=\delta_{ij},\\
&&\{\psi(i),\psi(j)\}=0=\{\psi^*(i),\psi^*(j)\},
\end{eqnarray*}
where we have put $\{a,b\}=ab+ba$.
We split $W$ into two subspaces given by
$W_a=(\oplus_{i<0}\Q\psi(i))\oplus(\oplus_{i\geq 0}\Q\psi^*(i))$
and
$W_c=(\oplus_{i\geq 0}\Q\psi(i))\oplus(\oplus_{i< 0}\Q\psi^*(i))$.
Then the Fock space (resp. the dual Fock space) is the
left (resp. right) $\Cl$-module $\F=\Cl/\Cl W_a=\Cl|0\>$
(resp. $\F^*=W_c \Cl\backslash \Cl=\<0|\Cl$), where
$|0\>=1$ mod $\Cl W_a$ (resp. $\<0|=1$ mod $W_c\Cl$).
A linear form $\<0|a|0\>\in \Q$ $(a\in \Cl)$ is uniquely determined
by setting $\<0|1|0\>=1$.

Define the charge zero sector $\F[0]$ as a subspace of $\F$
spanned by the elements
$\Psi=\psi^*(i_1)\psi^*(i_2)\cdots\psi^*(i_N)\psi(j_1)\psi(j_2)
\cdots\psi(j_N)|0\>$
$(j_N > \cdots > j_1\geq 0 > i_N > \cdots >i_1)$, namely
($\#$ of $\{\psi\})=(\#$ of $\{\psi^*\})$.
The space $\F[0]$ affords an irreducible representation of the Lie algebra
$\gl(\infty)=\gtg\oplus\Q\cdot 1$, where
$\gtg=\{\sum a_{ij}:\psi(i)\psi^*(j):\mid a_{ij}\in \Q,\,\,\,
a_{ij}=0 \mbox{ if }\\ |i-j|\gg 0\}$.
Note that the sum in the definition is infinite sum.

The boson-fermion correspondence gives an isomorphism
$\F[0]\simeq \Q[t_1,t_2,\ldots]$ as $\gl(\infty)$-modules.
More precisely, the basis vector above corresponds to,up to sign,
the Schur function associated with the generalized Maya diagram
$(j_N,\ldots,j_1,\bar{i}_1,\bar{i}_2,\ldots)$, where
$\{\bar{i}_k\}_{k\in \Z}=\{-1,-2,\ldots \}\backslash \{i_1,i_2,\ldots,i_N\}$
and $\bar{i}_1>\bar{i}_2>\cdots$.

As is shown in \cite{JM}, the Chevalley generators of $\A$
are realized in $\gl(\infty)$ as follows:
\begin{eqnarray*}
& & e_i=\sum_{j\in r\Z+i} \psi(j-1)\psi^*(j), \,\,\,\,\,
    f_i=\sum_{j\in r\Z+i} \psi(j)\psi^*(j-1), \\
& & \a_i^{\vee}=\sum_{j\in r\Z+i}
       \left(:\psi(j-1)\psi^*(j-1):-:\psi(j)\psi^*(j):\right)+\delta_{i0}
       \,\,\,(0\leq i\leq r-1),
\end{eqnarray*}
where : : denotes the normal ordering defined by
$$
:\psi(i)\psi^*(j):=\psi(i)\psi^*(j)-\<0|\psi(i)\psi^*(j)|0\>.
$$
The vacuum $|0\>$ generates an $\A$-module in $\F[0]$ by the action of the
Chevalley generators given above.
It is shown in \cite{JM} that this is the basic $\A$-module.
Applying the boson-fermion correspondence to this case, we see that
the representation space turns out to be $P^{(r)}$.

Let $\delta=\sum_{i=0}^{r-1} \a_i$ be the fundamental imaginary root.
If we denote by $W$ the Weyl group of type $\A$,
then the set of weights of the basic $\A$-module is given by
$$
P=\{w\L_0-n\delta \mid w\in W,\,n\in \N \}.
$$
It is known that the maximal weight vectors are the Schur functions
indexed by the $r$-core Young diagrams $\l_c$.
We denote by $\L(\l_c)$ the corresponding maximal weight.
Remark that these Schur functions coincide with the
$r$-reduced Schur functions.
Then we have the following theorem.

\begin{theorem}
The set of the $r$-reduced Schur functions
$$
\left\{S^{(r)}_{\l}(t) \mid \l=(\l_c,\vac,\l[1],\ldots,\l[r-1]),\,\,
\sum_{k=1}^{r-1} |\l[k]|=n\right\}
$$
gives a basis for a weight space of the basic $\A$-module with
a weight $\L(\l_c)-n\delta$.
\end{theorem}
\noindent
{\bf Proof.}
A key for the proof is that
to append a $r$-hook to $\l$ is nothing but to let \linebreak
$\psi(m+r)\psi^*(m)$ act on
$\Psi$ for some $m\in \Z$.

We prove that
$[\a_i^{\vee},\psi(m+r)\psi^*(m)]=0$ for $i=0,\ldots,r-1$ and $m\in \Z$.
However, this is a direct consequence of the identity
$$[\psi(k)\psi^*(k), \psi(m+r)\psi^*(m)]=
\psi(k)\psi^*(m)\delta_{k,m+r}-\psi(m+r)\psi^*(k)\delta_{mk}.
$$
On the other hand, $d_0$ acts on weighted homogeneous polynomials as
the degree counting operator.
Therefore, the $r$-reduced Schur function
$S^{(r)}_{\l}(t)$, $\l=(\l_c,\vac,\l[1],\ldots,\l[r-1])$,
is a weight vector of weight $\L(\l_c)-n\delta$ if
$\sum_{k=1}^{r-1} |\l[k]|=n$.
We can show that these elements span the weight space by
looking at the multiplicity formula
$$
\sum_{n=0}^{\infty} \mbox{mult}(\L(\l_c)-n\delta)q^n
=\frac{1}{\phi(q)^{r-1}}.
\mbox{\hspace{12pt}}\blacksquare
$$

\stepcounter{sect}
\begin{flushleft}
{\large\bf \S \thesect \hspace{5mm} A Formula for Modular Characters}
\end{flushleft}
\vspace{3mm}

The formula (Theorem 3.1) can be interpreted as a formula of the
Brauer characters of symmetric groups when $r$ is a prime number.
Let $(K,{\cal O},k)$ be an $r$-modular system.
In accordance
with the decomposition of the identity into the sum of primitive
central idempotents, the group algebra $k\Sym_n$ decomposes
into the direct sum of blocks, and the category of
$k\Sym_n$-modules decomposes into the
direct sum of the module categories of these blocks.
By the general theory on Specht modules $S^{\l}$,
it is well known \cite{JK,NT} that
\begin{enumerate}
\item Reduction modulo $r$ of any ${\cal O} \Sym_n$-submodule of an
irreducible $k \Sym_n$-module $S^{\l}$ ($|\l|=n$) has the same
set of composition factors all belonging to a single block.
\item Blocks are parametrized by $r$-cores and $S^{\l}$ belongs
to the block corresponding to the $r$-core $\l_c$ of $\l$.
\item The Grothendieck group of the module category of the block
corresponding to $\kappa$ is generated by
$\{[S^{\l}] \mid \l=(\l_c; \l[0],\ldots,\l[r-1]) \}$ as a
$\Z$-module.
\item Let $\chi^{\l}$ be the ordinary character of $S^{\l}$.
We consider it as a class function on $r$-regular classes.
Then the Grothendieck group is isomorphic to the $\Z$-submodule
of the space of class functions on $r$-regular classes spanned by
$\{\chi^{\l} \mid \l=(\l_c; \l[0],\ldots,\l[r-1]) \}$, in which
$[S^{\l}]$ corresponds to $\chi^{\l}$.
\end{enumerate}

Now we can state a theorem.
Because of the formula
$S_{\l}=\sum_{\mu} (1/z_{\mu}) \chi^{\l}(\mu) p_{\mu}$,
our formula (Theorem 3.1) is nothing but an equation for class functions
on $r$-regular classes. Hence we have

\begin{theorem}
\begin{enumerate}
\item $\{\chi^{(\l_c; \vac,\l[1],\ldots,\l[r-1])} \mid
|\l_c|+r\displaystyle{\sum_{k=1}^{r-1}| \l[k]|=n} \}$ is a
basic set.
\item The following equation holds on $r$-regular classes:
$$
\chi^{(r)}_{\l}=(-1)^{|\l[0]|} \delta_r(\l)
\! \! \! \sum_{\mu,\nu_1,\ldots,\nu_{r-1}} \! \! \!
LR^{\l[0]'}_{\nu_1\cdots \nu_{r-1}} LR^{\mu[1]}_{\nu_1 \l[1]} \cdots
LR^{\mu[r-1]}_{\nu_{r-1} \l[r-1]} \delta_r(\mu)
\chi^{(r)}_{\mu},
$$
where summation runs over Young diagram $\mu$ and $\nu_1,\ldots,\nu_{r-1}$
such that $|\mu|=|\l|$, $\mu[0]=\vac$ and the core of $\mu$ coincides
with that of $\l$.
\end{enumerate}
\end{theorem}

As is pointed out by A. Gyoja \cite{G}, we can lift this situation to
the Hecke algebra $H_n(q)$ with $q=\zeta_r$, a primitive $r$-th root of unity.
The $K$-algebra $H_n(\zeta_r)$ is defined by generators
$T_1, \ldots, T_{n-1}$ and relations
\begin{eqnarray*}
& &(T_i+1)(T_i-\zeta_r)=0 \,\,\,(1\leq i \leq n-1), \\
& & T_i T_{i+1} T_i=T_{i+1} T_i T_{i+1} \,\,\,(1\leq i \leq n-2), \\
& &T_i T_j=T_j T_i \,\,\,(|i-j|\geq 2).
\end{eqnarray*}
The reduction modulo $r$ of the ${\cal O}$-lattice
$\sum_{w\in \sym_n} {\cal O} T_w$ is isomorphic to $k \Sym_n$.
Hence we have a canonical $\Z$-linear map $\phi$ from the
Grothendieck group of $H_n(\zeta_r)$-modules to that of $k \Sym_n$-modules.
By \cite[Theorem 7.7]{DJ}, the decomposition matrices of both
are lower unitriangular with respect to the same ordering of
partitions. Hence we have
$$
\left[
      \begin{array}{ccc}
	1 &  & \bigzerou \\
	  & \ddots &     \\
	\bigast  &  & 1         \\
	  & \bigast &
      \end{array}
\right]
\left[
      \begin{array}{c}
	[D_{\l_1}] \\
	\hbox{$[D_{\l_2}]$} \\
	\vdots   \\
	\hbox{$[D_{\l_l}]$}
      \end{array}
\right]=
\left[
      \begin{array}{ccc}
	1 &  & \bigzerou \\
	  & \ddots &     \\
	\bigast &  & 1         \\
	  & \bigast &
      \end{array}
\right]
\left[
      \begin{array}{c}
	[E_{\l_1}] \\
	\hbox{$[E_{\l_2}]$} \\
	\vdots   \\
	\hbox{$[E_{\l_l}]$}
      \end{array}
\right],
$$
where $\{D_{\l_1},\ldots,D_{\l_l}\}$ are irreducible $k \Sym_n$-modules
and $\{E_{\l_1},\ldots,E_{\l_l}\}$ are irreducible $H_n(\zeta_r)$-modules.
By multiplying
$
\left[
\left.
      \begin{array}{ccc}
	1 &  & \bigzerou \\
	  & \ddots &     \\
	\bigzerol &  & 1         \\
      \end{array}
\right|
     \,	\hbox{\bg 0} \,
\right]
$
to the both sides, we have
$$
\left[
      \begin{array}{c}
	[D_{\l_1}] \\
	\vdots     \\
	\hbox{[$D_{\l_l}$]}
      \end{array}
\right]=
\left[
      \begin{array}{ccc}
	1 &  & \bigzerou \\
	  & \ddots &     \\
	\bigast &  & 1
      \end{array}
\right]
\left[
      \begin{array}{c}
	[E_{\l_1}] \\
	\vdots   \\
	\hbox{[$E_{\l_l}$]}
      \end{array}
\right].
$$
Therefore $\phi$ is invertible and hence these Grothendieck
groups are isomorphic.

\begin{corollary}
$$
[S^{(r)}_{\l}]=(-1)^{|\l[0]|} \delta_r(\l)
\! \! \! \sum_{\mu,\nu_1,\ldots,\nu_{r-1}} \! \! \!
LR^{\l[0]'}_{\nu_1\cdots \nu_{r-1}} LR^{\mu[1]}_{\nu_1 \l[1]} \cdots
LR^{\mu[r-1]}_{\nu_{r-1} \l[r-1]} \delta_r(\mu)
[S^{(r)}_{\mu}],
$$
where $S^{\l}$ is the Specht module for the Hecke algebra $H_n(\zeta_r)$.
\end{corollary}

\end{document}